\documentclass{amsart}

\usepackage{amsmath,amsfonts,mathrsfs,amssymb,amsthm,latexsym,array}
\usepackage[all]{xy}

\newcommand{\N}{\mathcal{H}} 
\newcommand{\A}{\mathrm{A}} 
\newcommand{\SI}{\mathrm{I}} 
\newcommand{\CM}{{\overline{\mathscr{M}}_{\N}}} 
\newcommand{\M}{{\mathscr{M}_{\N}}} 
\newcommand{\C}{\mathcal{C}} 
\newcommand{\D}{\mathcal{D}} 
\newcommand{\SQ}{\mathcal{Q}_{\N}} 
\newcommand{\X}{\mathcal{X}} 
\newcommand{\LO}{\mathcal{O}} 
\newcommand{\T}{\mathcal{T}} 
\newcommand{\G}{\mathbb{G}} 
\renewcommand{\k}{\mathbf{k}} 
\newcommand{\proj}{\mathbb{P}} 
\newcommand{\Nat}{\mathbb{N}} 
\renewcommand{\l}{\ell} 
\newcommand{\F}[1]{F_{#1}} 
\newcommand{\AF}{\mathbb{A}} 

\newcommand{\vv}{\mathfrak{o}} 
\renewcommand{\sharp}{\#}

\DeclareMathOperator*{\codim}{codim} 
\DeclareMathOperator*{\Spec}{Spec} %

\newtheorem{thm}{Theorem}[section]
\newtheorem{lem}[thm]{Lemma}
\newtheorem{cor}[thm]{Corollary}

\newtheorem{rmk}[thm]{Remark}
\newtheorem{syslem}[thm]{Syzygy Lemma}

\def \a[#1]{{g_{4,#1}}}
\def \b[#1]{{g_{5,#1}}}
\def \cc[#1]{{f_{6,#1}}}
\def \d[#1]{{g_{6,#1}}}
\def \e[#1]{{g_{7,#1}}}
\def \h[#1]{{f_{8,#1}}}
\def \kk[#1]{{g_{8,#1}}}
\def \m[#1]{{f_{9,#1}}}
\def \n[#1]{{f_{10,#1}}}

\def \f[#1,#2]{{f_{#1}^{(#2)}}}
\def \g[#1,#2]{{g_{#1}^{(#2)}}}

\def \ff[#1,#2]{{\tilde{f}_{#1}^{(#2)}}}
\def \gg[#1,#2]{{\tilde{g}_{#1}^{(#2)}}}

\begin{document}

\title[An upper bound for the dimension of $\M$]
{Upper bounds for the dimension of moduli spaces of
curves with symmetric Weierstrass semigroups}

\author{Andr\'e Contiero}
\address{Universidade Federal de Alagoas -- Instituto de Matem\'atica,
Avenida Lourival de Melo Mota, 57072-970 -- Macei\'o --AL, Brazil }

\email{andrecontiero@mat.ufal.br}

\thanks{The first author was supported by FAPEAL, Brazil}

\author{Karl-Otto St\"ohr}
\address{Instituto de Matem\'atica Pura e Aplicada (IMPA),
Estrada Dona Castorina 110, 22460-320 -- Rio de Janeiro -- RJ, Brazil}
\email{stohr@impa.br}


\subjclass[2010]{Primary 14H10, 14H20, 14H55}

\keywords{Moduli spaces, Weierstrass gaps, Gorenstein curves}

\date{ \today }


\begin{abstract}
We present an explicit 
method to produce upper bounds for the
dimension of the moduli spaces of complete integral Gorenstein
curves with prescribed symmetric Weierstrass semigroups.

\end{abstract}

\maketitle

\section{Introduction}
Let $\M$ be the moduli space of the smooth complete integral pointed
algebraic curves with a prescribed Weierstrass semigroup $\N$ of genus $g$.
There are two important estimates on the dimension of $\M$.
On the one hand,
Eisenbud and Harris \cite{EH87}, arguing that locally $\M$
is the pullback of Schubert-cycles from a suitable Grassmannian of $(g{-}1)$-planes,
obtained the lower bound $3g-2-w(\N)$ for the dimension of any irreducible
component of $\M$\,, where $w(\N)$ denotes the weight of $\N$\,.
As follows from their theory of limit linear series, this bound is
attained for some component of $\M$\,,
if $w(\N) \leq \frac{g}{2}\,,$ or more generally, if
$\N$ is a primitive semigroup of weight smaller than $g$  
(cf.\ \cite[theorem 3]{EH87}, \cite{K91}). However, if the weight is large,
as in the case of symmetric semigroups, then their bound may be 
far from being sharp, and it may even be negative.

On the other hand,
a theorem of Deligne \cite[theorem 2.27]{Del73}, whose proof
involves an interplay between 
three different moduli spaces, 
provides 
the upper bound $\dim \M \leq 2g-2+\lambda(\N)$, where
$\lambda(\N) \geq 1$ stands for the number of gaps $\l$ such that $\l+n$
is a nongap for each positive nongap $n$ (cf.\ \cite[section 6]{RV77}).
If the semigroup is symmetric, then Deligne's upper bound is equal to $2g-1$.
By using the work of Kontsevich and Zorich \cite{KZ03}, it has been
noticed by Bullock \cite{Bu13} that for each $g>3$ the upper bound
is attained exactly on the three symmetric semigroups with the gap
sequences $\{1,\dots,g-1,2g-1\}$, $\{1,\dots,g-2,g,2g-1\}$ and
$\{1,3,5,\dots,2g-1\}$.

In this paper we assume that $\N$ is a numerical symmetric semigroup.
Looking for an upper bound for the dimension of the moduli space
$\M\,,$ we view $\M$ as an open 
subspace of the compactified moduli space $\CM$ which is
defined by allowing arbitrary Gorenstein singularities.
By varying the construction of $\CM$ presented in \cite{St93},
in Section \ref{construction}
we realize $\CM$ in a rather explicit way as the weighted
projectivization $\proj(\X_{\N})$ of an affine quasi-cone $\X_{\N}$\,,
i.e.\ of a subset $\X_{\N}$ of a weighted vector space cut
out by quasi-homogeneous equations. 
We approximate the quasi-cone $\X_{\N}$ at its vertex
by a quadratic one, whose
weighted projectivization provides us with an upper bound for the
dimension of $\CM$ (cf.\ Theorem \ref{ubound}). We explain that it is
much easier to compute the quasi-homogeneous quadratic
forms that determine the quadratic approximation than the equations
of the moduli space $\CM\,$.


In the last section we illustrate the method by handling explicitely with
the family of symmetric semigroups
$$
  \N_\tau = \,<\! 6,2+6\tau,3+6\tau,4+6\tau,5+6\tau \!>
$$
of genus $\,g_\tau = 6\,\tau +1 \ (\tau =1,2,\dots)$.
The quadratic approximation of the quasi-cone ${\mathcal X}_{{\N}_\tau}$
can be described in terms of a variety defined over an artinean algebra
(cf.\ Theorem \ref{aq6}). Its dimension can be read off from this
description, providing the upper bound $\,8\,\tau +5\,$ for the
dimension of $\overline{\mathscr{M}}_{{\N}_\tau}\,,$
which for each $\tau > 1$ is better than Deligne's bound
$\,2\, g_\tau - 1 = 12\,\tau +1\,$.

\section{On the construction of the compactified moduli space
$\CM$}\label{construction}

Let $\C$ be a complete integral Gorenstein curve of arithmetic genus $g$ defined
over an algebraically closed field $\k$, and let $P$ be a smooth point of $\C$.
We denote by $\N$ the Weierstrass semigroup of the nongaps $0=n_0<n_1<n_2<\dots$
of the pointed curve $(\C,P)$. Thus for each $n\in\N$ there is a rational
function $x_n$ on $\C$ with pole divisor $nP$. We can assume that $x_0=1$.
For each nonnegative integer $i$ the vector space of global sections of the
divisor $n_i\,P$ is equal to
$$H^0(\C,n_i\,P)=\k\,x_{n_0}\oplus\k\,x_{n_1}\oplus\dots\oplus\k\,x_{n_i}$$
and in particular it has dimension $i+1$. It follows from Riemann's theorem for
complete integral
curves with singularities that $$n_{i}=g+i\ \ \text{for each }i\geq g$$ and so
the Weierstrass gap sequence consists of $g$ elements, say
$$1=\l_1<\l_2<\dots<\l_g\leq 2g-1.$$
We suppose that the last nongap $\l_g$ is equal to $2g-1$. This means that the
semigroup $\N$ is symmetric, or equivalently, it has the property that a
positive integer $\l$ is a gap if and only if $\l_g-\l$ is a nongap, i.e.
$$n_i=2g-1-\l_{g-i}\ \ (i=0,\dots,g-1)\,.$$
Thus $n_{g-1}=2g-2$ and $H^0(\C,(2g-2)P)$ is spanned by the $g$ functions
$x_{n_0}$, $x_{n_1},\dots,x_{n_{g-1}}$. Hence $\dim H^0(\C,(2g-2)P)=g$ and
so $(2g-2)P$ is a canonical divisor.
We also suppose that $\l_2=2$, or equivalently, the Weierstrass point $P$ is
nonhyperelliptic. Therefore by a theorem of Rosenlicht the canonical morphism
of the complete integral Gorenstein curve $\C$ is an embedding
$$(x_{n_0}:x_{n_1}:\dots:x_{n_{g-1}}):\C\hookrightarrow \proj^{g-1}$$
(see \cite[theorem 4.3]{KM09}). Thus $\C$ becomes a curve of degree $2g-2$ in
the projective space $\proj^{g-1}$ and the integers $\l_i-1$ $(i=1,\ldots,g)$
are the contact orders of the curve with the hyperplanes at $P=(0:\dots:0:1)$.

Conversely, the numerical symmetric nonhyperelliptic semigroup $\N$ can be
realized as the Weierstrass semigroup of the canonical monomial curve
$$\C^{(0)}:=\left\lbrace (a^{n_0}b^{\l_g-1}:a^{n_1}b^{\l_{g-1}-1}:\dots:
a^{n_{g-1}}b^{\l_1-1})\mid (a:b)\in\proj^{1}\right\rbrace\subset\proj^{g-1}$$
at the point $P$. The curve $\C^{(0)}$ is rational and has a unique singularity,
namely the unibranch point $(1:0:\dots:0)$ of multiplicity $n_1$ and singularity
degree $g$ (see \cite[p.\ 190]{St93}).

To study relations between generators of the ideal of the canonical curve
$\C\subset\proj^{g-1}$, we consider the spaces of global sections
$H^{0}(\C,\, r(2g-2)P)$ of the multi-canonical divisors $r(2g-2)P$, and
construct $P$-hermitian $r$-monomial bases, i.e. bases
consisting of $r$-monomial expressions in $x_{n_0},\dots,x_{n_{g-1}}$ whose
$P$-orders are pairwise different (cf.\ \cite{O91}).

\begin{lem}\label{basisbicanon}
Let $\tau$ denote the largest integer such that $n_{\tau}=\tau\,n_1$.
A $P$-hermitian basis of $H^{0}(\C,\, (4g-4)P)$ is given by the products
\begin{equation*}
\begin{array}{ll}
 x_{n_0}x_{n_j} & (j=0,\dots,g-1) \\
 x_{n_i}x_{n_{g-j}} & (i=1,\dots,\tau,\ \ j=1,\ldots,n_1-1) \\
 x_{n_i}x_{n_{g-j}} & (i=\tau+1,\dots,g-1,\ \ j=1,\dots,n_{i}-n_{i-1}) \\
 x_{n_{\tau+1}}x_{\l_{g-\tau-1}+k\,n_{1}} & (k=1,\dots,\tau-1) \, .
\end{array}
\end{equation*}
\end{lem}
\begin{proof}
Since by Riemann's theorem the dimension of $H^{0}(\C,\, (4g-4)P)$ is equal
to $3g-3$, we have only to convince ourselves that the above products are
$3g-3$ in number and have pairwise different pole orders at $P$
(see \cite[proposition 1.4]{St93}).
\end{proof}

For each nongap $s\leq 4g-4$ we list the partitions of $s$ as sums of two
nongaps $\leq 2g-2$, say
$$s=a_{si}+b_{si} \ \ \text{where}\ \ a_{si}\leq b_{si} \ \
(i=0,\dots,\nu_{s})$$
and abbreviate
$$a_s:=a_{s0} \ \ \text{and}\ \ b_s:=b_{s0}\,.$$
The $3g-3$ products $x_{a_s}x_{b_s}$ form a $P$-hermitian basis of the space of
global sections $H^0(\C,\,(4g-4)P)$ of the bi-canonical divisor $(4g-4)P$. If
$a_{s}<a_{s1}<\dots < a_{s\nu_{s}}$ then it coincides with the basis
displayed in Lemma \ref{basisbicanon}. However, we will not always work with
this assumption.

\begin{lem}\label{basismulticanon}
For each integer $r\geq 3$, a $P$-hermitian basis of $H^0(\C,\, r(2g-2)P)$ is
given by the $r$-monomial expressions
\begin{equation*}
\begin{array}{ll}
 x_{n_j} & (j=0,\dots,g-1) \\
 x_{a_s}x_{b_s}x^{i}_{n_{g-1}} & (i=0,\dots,r-2, \ \ s=2g,\dots,4g-4) \\
 x_{n_1}x_{2g-n_1}x_{n_{g-2}}x^i_{n_{g-1}} & (i=0,\dots,r-3)
\end{array}
\end{equation*} where the powers of $x_0=1$ have been omitted.
\end{lem}
\begin{proof}
Let $n$ be a nongap not larger then $r(2g-2)$. We look for a monomial
expression $z_n$ in $x_{n_1},\dots,x_{n_{g-1}}$ of degree $\leq r$ with pole
divisor $nP$. We proceed by induction on $r$. If $n\leq (r-1)(2g-2)$ then we
apply the induction hypothesis and pick up the corresponding basis element
of $H^0(\C,\,(r-1)(2g-2)P)$. Thus we may assume $(r-1)(2g-2)<n\leq r(2g-2)$.
If $n=4g-3$ then $r=3$ and we take $z_n:=x_{n_1}x_{2g-n_1}x_{2g-3}$. In the
remaining case we apply the induction hypothesis to $n-2g+2$ and take
$z_{n}:=z_{n-2g+2} \, x_{2g-2}$. \end{proof}

As becomes clear from the preceding proof,
we can normalize the functions $x_n$ in a way that for each $r\geq 2$ the functions
$x_n$ with $n\leq r(2g-2)$ are the above basis elements of $H^0(\C,r(2g-2)P)$.

Let $\SI(\C)\subset\k[X_{n_0},X_{n_1},\dots,X_{n_{g-1}}]$ be the homogeneous
ideal of the canonical curve $\C\subset \proj^{g-1}$ i.e.
\begin{equation*}
\SI=\bigoplus\limits_{r=2}^{\infty}\, \SI_{r}(\C)
\end{equation*}
where $ \SI_{r}(\C)$ is the vector space of the $r$-forms that vanish
identically on $\C$. As an immediate consequence of the existence of a
$P$-hermitian $r$-monomial basis for the space of global sections of the
multi-canonical divisor $r(2g-2)P$, the homomorphism
$$\k[X_{n_0},X_{n_1},\dots,X_{n_{g-1}}]_r\rightarrow H^0(\C,\, r(2g-2)P)$$
induced by the substitutions $X_{n_i}\mapsto x_{n_i}$ is surjective for each
$r$ (as predicted by a theorem of Noether).
Thus by Riemann's theorem the codimension of $\SI_{r}(\C)$ in the
$\binom{r+g-1}{r}$-dimensional vector space
$\k[X_{n_0},X_{n_1},\dots,X_{n_{g-1}}]_r$ of $r$-forms
is equal to $(2r-1)(g-1)$.
In particular, the vector space of quadratic relations has dimension
\begin{equation*}
\dim \SI_2(\C)=\frac{(g-2)(g-3)}{2} \, .
\end{equation*}
We attach to the variable $X_{n}$ the weight $n$. For each integer $r\geq 2$ we
fix a vector space $\Lambda_{r}$ in $\k[X_{n_0},X_{n_1},\dots,X_{n_{g-1}}]_r$
spanned by the lifting of a $P$-hermitian $r$-monomial basis of
$ H^0(\C,\, r(2g-2)P)$, or equivalently, a vector space spanned by
$r$-monomials in $X_{n_0},X_{n_1},\dots,X_{n_{g-1}}$ whose weights are
pairwise different and vary through the nongaps $n\leq r(2g-2)$. Since
$\ \Lambda_r\cap\SI_{r}(\C)=0$\, and
\begin{equation*}
\dim\Lambda_r=\dim  H^0(\C,\, r(2g-2)P)=\codim \SI_{r}(\C)\,,
\end{equation*}
we obtain
\begin{equation*}
\k[X_{n_0},X_{n_1},\dots,X_{n_{g-1}}]_r=\Lambda_r\oplus \SI_{r}(\C)
 \ \text{ for each } r\geq 2.
\end{equation*}
For each nongap $s\leq 4g-4$ and each integer $i=1,\dots,\nu_s$ we have
$x_{a_{si}}x_{b_{si}}\in H^0(\C,\, sP)$ and so we can write
\begin{equation*}
x_{a_{si}}x_{b_{si}}=\sideset{}{'}\sum_{n=0}^s\,c_{sin}x_{a_n}x_{b_n}
\end{equation*}
where the coefficients $c_{sin}$ are uniquely determined constants and where
the dash indicates that the summation index $n$ only varies through nongaps.
Multiplying the functions $x_{n_1},\dots,x_{n_{g-1}}$ by
constants we can normalize $$c_{sin}=1 \ \text{ whenever } \ n=s\,.$$
By construction the quadratic forms
\begin{equation*}
\F{si}:=X_{a_{si}}X_{b_{si}}-X_{a_s}X_{b_s}
  -\sideset{}{'}\sum_{n=0}^{s-1}\,c_{sin}X_{a_n}X_{b_n}
\end{equation*}
vanish identically on the canonical curve $\C$. They are linearly independent,
their number is equal to $\binom{g+1}{2}-(3g-3)=\frac{1}{2}(g-2)(g-3)$,
and hence they form a basis of the vector space of quadratic relations
$\SI_{2}(\C)$.

We attach to each coefficient $c_{sin}$ the weight $s-n$. If we consider $\F{si}$
as a polynomial expression not only in the variables $X_n$ but also in the
coefficients $c_{sin}$, then it becomes quasi-homogeneous of weight $s$.

To assure that the canonical ideal $\SI(\C)$ is generated by quadratic
relations, we need assumptions on the semigroup $\N$. In the following
we will always suppose that the symmetric semigroup $\N$ satisfies
$$3<n_1<g \ \ \text{ and } \ \ \N\neq\;<\!4,5\!>$$ i.e. $\N$ is different
from $\Nat\setminus\{1,\dots,g-1,2g-1\}$, $<\!2,2g+1\!>$, $<\!3, g+1\!>$ and $<\!4,5\!>$.
In the excluded case where $n_1=3$, respectively $\N=\;<\!4,5\!>$, the
canonical curve $\C\subset\proj^{g-1}$ is trigonal, respectively isomorphic to
a plane quintic, and one may easily check that the intersection of the
quadratic hypersurfaces containing $\C$ is an algebraic surface, or more
precisely, a rational normal scroll, respectively the Veronese surface
in $\proj^5$.

By a theorem of Oliveira \cite[theorem 1.7]{O91} it follows from the
assumptions on the
symmetric semigroup $\N$ that $\nu_{s}\geq 1$ whenever $s=n_i+2g-2$ and
$i=0,\dots,g-3$. Thus we can assume that
$$(a_{s1},b_{s1})=(n_i,2g-2) \ \
\text{whenever} \ \  s=n_i+2g-2 \ \ \text{and} \ \ i=0,\dots,g-3 \, .$$
Moreover, in this case we may assume that
$b_s=b_{s0}$ is the largest nongap smaller than $n_i+\l_{g-i-1}$ (see
\cite[proposition 1.7]{St93}).

The canonical curve $\C\subset \proj^{g-1}$ is contained in the $g-2$ quadratic
hypersurfaces $V(\F{n_i+2g-2,1})$ $(i=0,\dots,g-3)$, which intersect
transversely
at the Weierstrass point $P$. Thus in an open neighborhood of $P$, the
curve $\C$
is the intersection of these $g-2$ hypersurfaces. Therefore the $g-2$ quadratic
forms $\F{n_i+2g-2,1}$ determine uniquely the integral curve $\C$, and in
particular determine the remaining $\frac{1}{2}(g-2)(g-5)$ quadratic forms of
the basis of $\SI_{2}(\C)$. We will make this explicit, by
constructing syzygies
of the canonical curve.

It now follows from the Enriques-Babbage theorem that the canonical integral
curve $\C$ is non-trigonal and not isomorphic to a plane quintic (see \cite[p.\
124]{ACGH85}). Moreover by Petri's analysis, at least in the case where $\C$ is
smooth, the canonical ideal $\SI(\C)$ is generated by quadratic relations
(see \cite[p.\ 131]{ACGH85}).

We will give an algorithmic proof that the ideal $\SI(\C)$ of the (possibly
nonsmooth) canonical curve $\C\subset\proj^{g-1}$ is generated by the quadratic
forms $\F{si}$. We first treat the canonical monomial curve $\C^{(0)}$.

\begin{lem}\label{ic0}
The ideal $\SI(\C^{(0)})$ of the canonical monomial curve
$\C^{(0)}\subset\proj^{g-1}$ is generated by the quadratic binomials
$$\F{si}^{(0)}:=X_{a_{si}}X_{b_{si}}-X_{a_s}X_{b_s}\, .$$
\end{lem}
\begin{proof}
Let $\SI^{(0)}$ be the ideal generated by the binomials $\F{si}^{(0)}$. Since
the ideal $\SI(\C^{(0)})$ is homogeneous and quasi-homogeneous, it is enough to
show that a homogeneous and quasi-homogeneous polynomial belongs to $\SI^{(0)}$
if it belongs to $\SI(\C^{(0)}),$ i.e.\ if the sum of its coefficients is equal to zero.

In this proof we do not need the assumption $a_{n_i+2g-2,1}=n_i \, ,$ and so without
loss of generality we can assume that $a_s<a_{s1}<\dots<a_{s\nu_s}$ for each
$s$.
We further assume that the vector spaces $\Lambda_r$ are obtained by lifting
the bases displayed in Lemma \ref{basisbicanon} and Lemma \ref{basismulticanon}.
To apply Gr\"obner basis techniques,  
we order the monomials $\prod_{k=0}^{g-1}X_{n_k}^{i_k}$ according to the
lexicographic ordering of the vectors
$(\sum\,i_k,\sum\,n_k\,i_k,-i_0,-i_{g-1},\dots,-i_1)$. Using our assumptions on
the semigroup $\N$\,, we can enlarge the basis $\{\F{si}^{(0)}\}$ of $\SI^{(0)}$
to a Gr\"obner basis 
by adding
certain sums of cubic binomials $\pm X_{n}\F{si}^{(0)}$ of the same weight
(see \cite[pp.\ 196--198]{St93}).

Let $F$ be a homogeneous polynomial of degree $r$. Dividing $F$ by the
Gr\"obner basis, the division algorithm provides us with a decomposition
$$F=\sum_{si}\,G_{si}\F{si}^{(0)}+R$$
where $R\in\Lambda_{r}$ and
$G_{si}$ is homogeneous of degree $r-2$ for each double index $si\,.$
If $F$ is quasi-homogeneous of weight $w$, then
each $G_{si}$ is quasi-homogeneous of weight $w-s$, and
the remainder $R$ is the only monomial in $\Lambda_{r}$ of weight
$w$ whose coefficient is equal to the sum of the coefficients of $F$.
Therefore,
if  $F\in\SI(\C^{(0)})$ then $R=0$ and $F\in\SI^{(0)}$.
\end{proof}

\begin{syslem}\label{syzlem}
For each of the $\frac{1}{2}(g-2)(g-5)$ quadratic binomials $\F{s'i'}^{(0)}$
different from 
$\F{n_i+2g-2,1}^{(0)}$ $(i=0,\dots,g-3)$ there is
a syzygy of the form
$$X_{2g-2}\F{s'i'}^{(0)}+\sum_{nsi}\varepsilon_{nsi}^{(s'i')}X_{n}\F{si}^{(0)}=0$$
where the coefficients $\varepsilon_{nsi}^{(s'i')}$ are integers equal to $1$,
$-1$ or $0$ (which will be specified below), and where the sum is
taken over the
nongaps $n<2g-2$ and the double indices $si$
with $n+s=2g-2+s'$.
\end{syslem}

A weak version of the Syzygy Lemma 
has been obtained in \cite[lemma
2.3]{St93} by using Petri's analysis. We will provide an elementary
purely combinatorial proof.

\begin{proof}
If we put $F:=\F{s'i'}^{(0)}$ or  $F:=-\F{s'i'}^{(0)}$ then we can write
$$F=X_{q}X_{r}-X_{m}X_{n}$$
where $q$, $r$, $m$ and $n$ are nongaps satisfying $q+r=m+n$ and
$m<q\leq r < n < 2g-2\,.$  The strict inequality $n<2g-2$ is due to the
assumption $b_{s1}=2g-2$ whenever $s=n_i+2g-2$.

If $n+1$ is a gap then, by symmetry, the integer $k:=2g-2-n+r$ is a nongap
smaller than $2g-2$, and we have the syzygy
$$X_{2g-2}F+X_{n}(X_{m}X_{2g-2}-X_{q}X_{k})-X_{q}(X_{r}X_{2g-2}-X_{n}X_{k})=0$$
where the binomials in the brackets can be written as
$\F{si}^{(0)}-\F{sj}^{(0)}$,
$\F{si}^{(0)}$ or $-\F{sj}^{(0)}$.
Analogously, if $q+1$ is a gap, then the integer $k:=2g-2-q+m$ is a nongap
smaller than $2g-2$, and we get the syzygy obtained from the previous one by
interchanging $n$ with $q$ and $r$ with $m$.
Now we can assume that $n+1$ and $q+1$ are nongaps. Then we have the syzygy

\medskip
\hspace{1cm}
$X_{2g-2}F + X_m(X_nX_{2g-2}-X_{n+1}X_{2g-3})=$

\smallskip
\hspace{1cm}
$X_{2g-3}(X_{q+1}X_{r}-X_mX_{n+1}) 
+ X_{r}(X_{q}X_{2g-2}-X_{q+1}X_{2g-3}) \, .$
\end{proof}

\smallskip

In the remainder of this section we invert the above considerations.
Given a numerical symmetric semigroup
$\N$ of genus $g:=\sharp(\Nat\setminus\N)$ satisfying $3<n_1<g$ and
$\N\neq\;<\!4,\!5>$, and given $\frac{1}{2}(g-2)(g-3)$ quadratic forms
$$\F{si}=\F{si}^{(0)}-\sideset{}{'}\sum_{n=0}^{s-1}\,c_{sin}X_{a_n}X_{b_n}\ ,$$
we look for conditions on the constants $c_{sin}$ in order that the
intersection of the quadratic hypersurfaces $V(\F{si})$ in $\proj^{g-1}$ is a
canonical integral Gorenstein curve.

\begin{lem}\label{ir+dr}
Let $\SI$ be the ideal generated by the $\frac{1}{2}(g-2)(g-3)$ quadratic
forms $\F{si}$. Then
$$\k[X_{n_0},X_{n_1},\dots,X_{n_{g-1}}]_r=\SI_{r}+\Lambda_{r} \ \
\text{for each} \ r\geq 2\ .$$
\end{lem}
\begin{proof}
Let $F$ be a homogeneous polynomial of degree $r$ and weight $w$. Let
$G$ be its
quasi-homogeneous component of weight $w$, and let $R$ be the only monomial in
$\Lambda_r$ of weight $w$ whose coefficient is equal to the sum of the
coefficients of $G$. Then $G-R\in\SI(C^{(0)})$, and so by Lemma \ref{ic0} there
is a partition
$$G-R=\sum_{si}G_{si}\F{si}^{(0)} \, .$$
Replacing each polynomial $G_{si}$ by its homogeneous component of
degree $r-2$,
we can assume that $G_{si}$ is homogeneous of degree $r-2$. For
similar reasons,
we can assume that each $G_{si}$ is quasi-homogeneous of weight $w-s$. Hence
the polynomial $F-\sum G_{si}\F{si}-R$ is homogeneous of degree $r$ and weight
$\leq w-1$. Now the lemma follows by induction on $w$. \end{proof}

Replacing in the left hand side of the Syzygy Lemma the binomials
$\F{s'i'}^{(0)}$ and $\F{si}^{(0)}$ by the quadratic forms $\F{s'i'}$ and
$\F{si}$ we obtain for each of the $\frac{1}{2}(g-2)(g-5)$ double
indices $s'i'$ a
linear combination of cubic monomials of weight $<s'+2g-2$, which by Lemma
\ref{ir+dr}, or more precisely, by its algorithmic proof admits a decomposition
\begin{equation*}\label{Rsi}
 X_{2g-2}\F{s'i'}+\sum_{nsi}\varepsilon_{nsi}^{(s'i')}X_{n}\F{si}=
 \sum_{nsi}\eta_{nsi}^{(s'i')}X_{n}\F{si}+R_{s'i'}
\end{equation*}
where the sum on the right hand side is taken over the nongaps $n\leq 2g-2$ and
the double indices $si$ with $n+s<s'+2g-2$, where the coefficients
$\eta_{nsi}^{(s'i')}$ are constants, and where $R_{s'i'}$ is a linear
combination of cubic monomials of pairwise different weights $<s'+2g-2$.

For each nongap $m<s'+2g-2$ we denote by $\varrho_{s'i'm}$ the only coefficient
of $R_{s'i'}$ of weight $m$. Since we do not loose information about the
coefficients of $R_{s'i'}$ if we replace the variables $X_{n}$ by powers $t^n$
of an indeterminate $t$, it is convenient to consider the polynomial
\begin{equation*}
R_{s'i'}(t^{n_0},t^{n_1},\dots,t^{n_{g-1}})=
\sideset{}{'}\sum_{m=0}^{s'+2g-3}\varrho_{s'i'm}t^m
\end{equation*}
where the dash indicates that the summation index only varies through nongaps.

Since $R_{s'i'}$ may be obtained as a remainder in a division procedure, we can
arrange that the coefficients $\varrho_{s'i'm}$ (respectively,
$\eta_{nsi}^{(s'i')})$ are quasi-homogeneous polynomial expressions of weight
$s'+2g-2-m$ (respectively, $s'+2g-2-n-s$) in the constants $c_{sin}$. However,
we do not specify a division procedure and we do not postulate that $R_{s'i'}$
belongs to $\Lambda_3$, and so the $\frac{1}{2}(g-2)(g-5)$ remainders $R_{s'i'}$
are not uniquely determinate. In practice, this freedom in the construction of
$R_{s'i'}$ allows us to make shortcuts,
as we will illustrate in the last section.

\begin{thm}
Let $\N\subset\Nat$ be a numerical symmetric semigroup of genus $g$ satisfying
$3<n_1<g$ and $\N\neq<4,5>$. Then the $\frac{1}{2}(g-2)(g-3)$ quadratic forms
$\F{si}=\F{si}^{(0)}-\sum_{n=0}^{s-1}c_{sin}X_{a_n}X_{b_n}$ cut out a canonical
integral Gorenstein curve  in $\proj^{g-1}$ if and only if the coefficients
$c_{sin}$ satisfy the quasi-homogeneous equations $\varrho_{s'i'm}=0$. In this
case, the point $P=(0:\dots:0:1)$ is a smooth point of the canonical curve with
Weierstrass semigroup $\N$.
\end{thm}
\begin{proof}
We first assume that the $\frac{1}{2}(g-2)(g-3)$ quadratic forms $\F{si}$ cut out
a canonical integral Gorenstein curve in $\proj^{g-1}$. Since each $R_{s'i'}$
belongs to the ideal $\SI$ generated by the quadratic forms $\F{si}$\,, we
have $R_{s'i'}(x_{n_0},x_{n_1},\dots,x_{n_{g-1}})=0$ for each double index
$s'i'$. On the other hand
$$R_{s'i'}(x_{n_0},\dots,x_{n_{g-1}})=
\sideset{}{'}\sum_{m=0}^{s'+2g-3}\varrho_{s'i'm}z_{s'i'm}$$
where each $z_{s'i'm}$ is a monomial expression of weight $m$ in the projective
coordinate functions $x_{n_0},\dots,x_{n_{g-1}}$, and hence has the
pole divisor $mP$. Thus the coefficients $\varrho_{s'i'm}$ are zero.

Let us now assume that the coefficients $c_{sin}$ satisfy the equations
$\varrho_{s'i'm}=0$. Since the $g-2$ quadratic hypersurfaces
$V(\F{n_i+2g-2,1})\subset\proj^{g-1}$ ($i=0,\dots,g-3$) intersect transversely
at $P$, their intersection has a unique irreducible component that passes
through $P$, and this component is a projective integral algebraic curve, say
$\C$, which is smooth at $P$ and whose tangent line at $P$ is the intersection
of their tangent hyperplanes $V(X_{n_i})$ ($i=0,\dots,g-3$).

Let $y_{n_0},\dots,y_{n_{g-1}}$ be the projective coordinate functions of $\C$
with $y_{n_{g-1}}=1$.
Since $n_{g-1} - n_{g-2} = \l_2 - \l_1 = 1,$ we conclude that  
$t:=y_{n_{g-2}}$ is a local parameter of $\C$ at $P$, and
$y_{n_o},\dots,y_{n_{g-3}}$ are power series in $t$ of order $>1$. More
precisely, comparing coefficients in the $g-2$ equations
$\F{n_i+2g-2,1}(y_{n_0},\dots,y_{n_{g-1}})=0$ we obtain
$$y_{n_i}=t^{n_{g-1}-n_i}+\dots\ =t^{\l_{g-i}-1}+\dots \ \ \ (i=0,\dots,g-1)$$
where the dots stand for terms of higher orders.
Thus the $g$ integers $\l_i-1$ $(i=1,\dots,g)$ are the contact orders of the
curve $\C\subset\proj^{g-1}$ with the hyperplanes. In particular, the curve
$\C$ is not contained in any hyperplane.

Since by assumption the $\frac{1}{2}(g-2)(g-5)$ remainders $R_{s'i'}$ are
equal to zero, we have the syzygies
$$X_{2g-2}\F{s'i'}+\sum_{nsi}\varepsilon_{nsi}^{(s'i')}X_{n}\F{si}-
\sum\eta_{nsi}^{(s'i')}X_{n}\F{si}=0\ .$$
Replacing the variables $X_{n_0},\dots,X_{n_{g-1}}$ by the projective
coordinate
functions $y_{n_0}$, $\dots$, $y_{n_{g-1}}$, we get a system of
$\frac{1}{2}(g-2)(g-5)$ linear homogeneous equations in the
$\frac{1}{2}(g-2)(g-5)$
functions $\F{s'i'}(y_{n_0},\dots,y_{n_{g-1}})$ with coefficients in the domain
$\k[[t]]$ of formal power series. Since the triple indices $nsi$ of the
coefficients $\varepsilon_{nsi}^{(s'i')}$, respectively $\eta_{nsi}^{(s'i')}$,
satisfy $n<2g-2$ and $n+s=s'+2g-2$, respectively $n\leq 2g-2$ and
$n+s<s'+2g-2$,
we conclude that the diagonal entries of the matrix of the system have
constant terms equal to $1$, while the remaining entries have
positive orders.
Thus the matrix is invertible, and so the equation
$\F{si}(y_{n_0},\dots,y_{n_{g-1}})=0$
holds for each double index $si$. This shows that $\SI\subseteq\SI(\C)$ where
$\SI$ is the ideal generated by the $\frac{1}{2}(g-2)(g-3)$ quadratic forms
$\F{si}\,$.

By Lemma \ref{ir+dr} we have $\codim\SI_{r}\leq \dim\Lambda_r$ for each
$r\geq 2$. On the order hand, since $\SI_r(\C)\cap\Lambda_r=0$, we deduce
$\dim\Lambda_r \leq \codim\SI_r(\C)$. 
Since $\SI\subseteq\SI(\C)$ we obtain
$$\codim\SI_r(\C)=\codim\SI_r=\dim\Lambda_r=(2g-2)r+1-g\ .$$
Thus $\SI(\C)=\SI$ and the curve $\C\subset \proj^{g-1}$ has Hilbert polynomial
$(2g-2)r+1-g$. Hence $\C$ has degree $2g-2$ and arithmetic genus equal to $g$.

Intersecting the curve $\C\subset\proj^{g-1}$ with the hyperplane
$V(X_{n_{g-1}})$ we get the divisor $D:=(2g-2)P$ of degree $2g-2$, whose
complete linear system $|D|$ has dimension $g-1$. Hence by the Riemann-Roch theorem
for complete integral (not necessarily smooth) curves the Cartier
divisor $D$ is canonical, and $\C$ is a canonical Gorenstein curve.
\end{proof}

To normalize the coefficients $c_{sni}$ of the quadratic forms $\F{si}$, we
notice that the coordinate functions $x_n$ ($n\in\N,\ n\leq 2g-2$) are not
uniquely determined by their pole divisors $nP$. We assume that the
characteristic of the constant field $\k$ is zero or a prime not
dividing any of
the differences $n-m$, where $m$ and $n$ are nongaps satisfying $m<n\leq 2g-2$.
Transforming
$$X_{n}\mapsto X_{n}+\sum_{m=0}^{n-1}d_{nm}X_{m}$$
where the coefficients $d_{nm}$ are constants, we can normalize
$\frac{1}{2}g(g-1)$ of the coefficients $c_{sin}$ to be zero. More precisely,
for each positive integer $w$
the number of the coefficients $c_{sin}$ of weight $s-n=w$ that can be
normalized is equal to the number of nongaps $m$ such that $m+w$ is a
nongap $\leq 2g-2$ (see \cite[proposition 3.1]{St93}).

Due to these normalizations and the normalizations of the coefficients of weight
zero (i.e.\ $c_{sis}=1$), the only freedom left to us is to transform
$x_{n_i}\mapsto c^{n_i}x_{n_i}$ ($i=0,\dots,g-1$) for some
$c\in\G_{m}(\k)=\k^{*}$.
We have shown:

\begin{thm}\label{moduli}
Let $\N\subset\Nat$ be a symmetric semigroup of genus
$g:=\sharp(\Nat\setminus\N)$
sa\-tis\-fying $3<n_1<g$ and $\N\;\neq<\!4,5\!>$. The isomorphism classes of the
pointed complete integral Gorenstein curves with Weierstrass semigroup $\N$
correspond bijectively to the orbits of the $\G_{m}(\k)$-action
$$(c,\dots,c_{sin},\dots)\mapsto (\dots,c^{s-n}c_{sin},\dots)$$
on the affine quasi-cone of the vectors whose coordinates are the coefficients
$c_{sin}$ of the normalized quadratic forms $\F{si}$ that satisfy the
quasi-homogeneous equations $\varrho_{s'i'm}=0$.
\end{thm}

\section{The Method}

We start this section by describing a variant of the construction of the
remainders $R_{s'i'}$. Instead of making induction on the weights in the
variables $X_n$, we proceed by induction on the degrees in the constants
$c_{sin}$. We choose for each pair $s'i'$ and each nongap $m\leq 6g-6$ a
cubic monic monomial $Z_{s'i'm}$ of weight $m$ 
in $X_{n_0},\dots,X_{n_{g-1}}$ e.g.\
the unique monic monomial in $\Lambda_{3}$ of weight $m$.
By the Syzygy Lemma the polynomial
\begin{equation*}
G_{s'i'}^{(1)}:=
X_{2g-2}\F{s'i'}+\sum_{nsi}\varepsilon_{nsi}^{(s'i')}X_n\F{si}\ ,
\end{equation*}
is a sum of cubic monomials of weight
$<s'+2g-2$, whose coefficients are homogeneous of degree $1$ in the constants
$c_{sin}$. For each nongap $m<s'+2g-2$ let $\varrho_{s'i'm}^{(1)}$ be
the sum of
the coefficients of its terms of weight $m$. We notice that
$\varrho_{s'i'm}^{(1)}$ is homogeneous of degree $1$ and quasi-homogeneous of
weight $s'+2g-2-m$ in the constants $c_{sin}$. By Lemma \ref{ic0}, or more
precisely, by its algorithmic proof, we get an equation
\begin{equation*}
G_{s'i'}^{(1)}-\sum_{m}\varrho_{s'i'm}^{(1)}Z_{s'i'm}=
\sum_{sim}\eta_{sim}^{(s'i')(1)}X_{m}\F{si}^{(0)}
\end{equation*}
where the coefficients $\eta_{sim}^{(s'i')(1)}$ are homogeneous of degree $1$
and quasi-homogeneous of weight $s'+2g-2-m-s>0$ in the constants
$c_{sin}$. Next
we consider the polynomial
\begin{equation*}
G_{s'i'}^{(2)}:=G_{s'i'}^{(1)}-\sum_{m}\varrho_{s'i'm}^{(1)}Z_{s'i'm} -
\sum_{sim}\eta_{sim}^{(s'i')(1)}X_{m}\F{si} =
\sum_{sim}\eta_{sim}^{(s'i')(1)}X_{m}(\F{si}^{(0)}-\F{si})\ ,
\end{equation*}
which is a sum of cubic monomials of weight $<s'+2g-2$, whose coefficients are
homogeneous of degree $2$ in the constants $c_{sin}$. Let
$\varrho_{s'i'm}^{(2)}$ be the sum of the coefficients of its terms of weight
$m$. Proceeding in this way, we obtain the coefficients of the remainders
$R_{s'i'}$ as sums of its homogeneous components:
\begin{equation*}
\varrho_{s'i'm}=\varrho_{s'i'm}^{(1)}+\varrho_{s'i'm}^{(2)}+\dots
\end{equation*}
The procedure stops after less than $s'+2g-2$ steps, because the degrees in the
constants are not larger than $s'+2g-2-m$.

To simplify the notation, we put the coefficients $c_{sin}$ of the normalized
quadratic forms $\F{si}$ in any order, rewrite them as a finite sequence
$c_0,c_1,\dots,c_{\nu}$ and denote their weights by $w_0,w_1,\dots,
w_{\nu}$. We also rewrite the equations $\varrho_{s'i'm}=0$ by using
simple indices
\begin{equation*}
\varrho_{j}(c_0,\dots,c_{\nu})=0 \ \ \ (j=1,\dots,\mu)
\end{equation*}
where each $\varrho_j$ is a quasi-homogeneous polynomial in $\nu+1$ variables,
whose homogeneous components $\varrho_{j}^{(1)}, \varrho_{j}^{(2)},\dots$ have
been determined by the above procedure. There is an action of the multiplicative
group $\G_{m}:=\k^{*}$ on the weighted vector space $\k^{\nu+1}$
defined as follows
\begin{equation*}
(c_0,\dots,c_{\nu})\mapsto (c^{w_0}c_0,\dots,c^{w_{\nu}}c_{\nu})
\ \ \text{for each} \ c\in\k^{*}\ .
\end{equation*}
By Theorem \ref{moduli} the isomorphism classes of the pointed complete integral
Gorenstein curves with Weierstrass semigroup $\N$ correspond bijectively
to the orbits of the $\G_{m}$-action on the affine quasi-cone
\begin{equation*}
\{(c_0,\dots,c_{\nu})\in\k^{\nu+1})\ |\ \varrho_{j}(c_0,\dots,c_{\nu})=0\ , \
j=1,\dots,\mu)\}\ .
\end{equation*}
Linearizing the $\mu$ equations we obtain a vector space
\begin{equation*}
\T_{\N}:=\{(c_0,\dots,c_{\nu})\in\k^{\nu+1}\ | \ \varrho_{j}^{(1)}=0 \
, \ j=1,\dots,\mu\}
\end{equation*}
equipped with the induced $\G_m$-action. To relate this
with Pinkham's work \cite{Pi74}, we notice that
\begin{equation*}
\T_{\N}\cong T^{1,-}_{\k[\N]|\k} \ ,
\end{equation*}
that is, $\T_{\N}$ is isomorphic to the negatively graded part of the
first coho\-mo\-logy group $T^{1}_{\k[\N]|\k}$ of the cotangent complex of the
semigroup algebra $\k[\N]|\k$, where the signs of the weights have been
inverted (see \cite[p.\ 212]{St93}). To determine the vector space $\T_{\N}$,
we eliminate successively coefficients $c_i$ from the linear equation
$\varrho_{j}^{(1)}(c_0,\dots,c_{\nu})=0$. Let $r$ be the smallest integer
such that, after an eventual  permutation, the coefficients $c_i$ with
$i>r$ can be eliminated from the linear equations. Then $\T_{\N}$ becomes
the $(r+1)$-dimensional weighted vector space
\begin{equation*}
\T_{\N}=\{(c_0,\dots,c_{r})\ |\ c_0,\dots,c_r\in\k\}
\end{equation*}
with the weight sequence $w_0,\dots,w_r$. If a coefficient $c_i$ with $i>r$ has
been eliminated from the linear equations
$\varrho_{j}^{(1)}(c_0,\dots,c_{\nu})=0$, then it can also be eliminated from the
corresponding polynomial equation $\varrho_{j}(c_0,\dots,c_{\nu})=0$, because by the
quasi-homogeneity it does not occur in higher order terms of the same equation.
Entering with these solutions into the remaining polynomials, after $\nu-r$
steps we obtain $\mu-(\nu-r)$ quasi-homogeneous polynomials $h_j$ in only $r+1$
variables, which define an affine quasi-cone
\begin{equation*}
\X_{\N}:=\{(c_0,\dots,c_r)\in\k^{r+1}\ | \ h_j(c_0,\dots,c_r)=0\ \forall j \}
 \subseteq \T_{\N}\ .
\end{equation*}
By Theorem \ref{moduli} the quotient $\X_{\N}\diagup\G_m$
parametrizes the isomorphism classes of the pointed complete integral Gorenstein
curves with Weierstrass semigroup $\N$.
Thus the compactified moduli space $\CM$ can be identified with the quotient
of the punctured quasi-cone $\X_{\N}\setminus\{\vv\}$ by the
$\G_{m}$-action
\begin{equation*}
\CM\cong\proj(\X_{\N}):=(\X_{\N}\setminus\{\vv\})
\diagup\G_m\subseteq\proj(\T_{\N})
\end{equation*}
i.e. $\CM$ is isomorphic to the closed subset of the $r$-dimensional
weighted projective space $\proj(\T_{\N})=\proj_{(w_0,\dots,w_r)}^{r}$
cut out by the quasi-homogeneous equations $h_j=0$. Here the pointed monomial
curve, which corresponds to the vertex $\vv$ of the affine quasi-cone
$\X_{\N}$, has been excluded, and can be viewed as the "improper point" of
$\CM$.

By construction, the linear components $h_{j}^{(1)}$ of the quasi-homogeneous
polynomials $h_j$ are equal to zero. The quadratic components $h_j^{(2)}$ can
be easily computed, by solving the linear equations $\varrho_{j}^{(1)}=0$ for
$c_{r+1},\dots,c_{\nu}$, and entering into the quadratic expressions
$\varrho_{j}^{(2)}(c_0,\dots,c_{r},c_{r+1},\dots,c_{\nu})$. We approximate
$\X_{\N}$ at the vertex by the affine quadratic quasi-cone
\begin{equation*}
\SQ:=\{(c_0,\dots,c_r)\in\k^{r+1}\ | \ h_{j}^{(2)}(c_0,\dots,c_r)=0
   \ \ \forall j\}\subseteq \T_{\N}\,.
\end{equation*}

\begin{thm}\label{ubound}
$$\dim\CM\leq \dim\proj(\SQ)\ \ \text{i.e.} \ \ \dim\CM<\dim\SQ\,.$$
\end{thm}
\begin{proof}
Since $\CM\cong(\X_{\N}\setminus\{\vv\})\diagup\G_{m}$ we have
$$\dim\CM=\dim\X_{\N}-1 \, .$$
Due to quasi-homogeneity each irreducible component of the affine quasi-cone
$\X_{\N}$ passes through the vertex $\vv$. Since the dimension of an
integral variety coincides with its local dimension at any point
(see \cite[\S 13, Theorem A]{Ein95}), we conclude that the dimension
of $\X_{\N}$
is equal to its local dimension at the vertex
$$\dim\X_{\N}=\dim_{\vv}\X_{\N}$$
which is equal to the Krull dimension of the corresponding local ring
$$\dim_{\vv}\X_{\N}=\dim\LO_{\X_{\N},\vv}\ .$$
Since by local algebra the dimension of a local ring is equal to the dimension
of its associated algebra (\cite[theorem 13.9]{M89}) we have
$$\dim\LO_{\X_{\N},\vv}=\dim G (\LO_{\X_{\N},\vv})$$
where
$$G(\LO_{\X_{\N},\vv}):=
\oplus_{i=1}^{\infty}\, ({\mathfrak m}_{\X_{\N},\vv})^i\diagup
({\mathfrak m}_{\X_{\N},\vv})^{i+1}\ .$$
Geometrically, this means that the local dimension is equal to the dimension
of the tangent cone
$$\dim_{\vv}\X_{\N}=\dim \mathrm{C}_{\vv}(\X_{\N})$$
and can be seen by noticing that the projectivization of the tangent cone is an
effective Cartier divisor in the local blowup. Since by construction the
quadratic quasi-cone $\SQ$ contains the tangent cone
$\mathrm{C}_{\vv}(\X_{\N})$, we conclude

\medskip
\hspace{1.5in}
$\dim\mathrm{C}_{\vv}(\X_{\N})\leq \dim\SQ\ .$
\end{proof}
It is much less expensive to obtain the equations and the dimension of the
quadratic quasi-cone $\SQ$ than the ones of the moduli space $\M$. With
Theorem \ref{ubound} we get an implementable method to produce an upper bound
for the dimension of the moduli space of curves with a prescribed symmetric
Weierstrass semigroup. Below we summarize in a table some examples we
calculated on a computer:

\begin{table}[ht!]\label{tablesample}
\centering
\begin{tabular}{|c|c|c|c|c|c|c|} \hline
$\N$ & $g$ & E-H & $\dim\CM$ & $\dim\proj(\SQ)$ & Del &
$\dim T^{1,-}_{\k[\N]|\k}$ \\ \hline
$<6,8,9,10,11>$ & $7$ & $12$ & $13$ & $13$ & $13$ & $17$ \\ \hline
$<6,8,10,11,13>$ & $8$ & $12$ & $14$ & $14$ & $15$ & $18$ \\ \hline
$<7,9,10,11,12,13>$ & $8$ & $14$ & $15$ & $15$ & $15$ & $23$ \\ \hline
$<6,8,10,13,15>$ & $9$ & $11$ & $15$ & $15$ & $17$ & $19$ \\ \hline
$<6,9,10,13,14>$ & $9$ & $12$ & $15$ & $15$ & $17$ & $19$ \\ \hline
$<6,14,15,16,17>$ & $13$ & $11$ & $20$ & $21$ & $25$ & $28$ \\ \hline
\end{tabular}
\end{table}
\noindent Here E-H stands for the lower bound $3g-2-w(\N)$ of
Eisenbud-Harris and Del for Deligne's upper bound $2g-1$.

By the jacobian criterion and classical elimination theory, the
moduli space $\M$ is an open subspace of
$\CM$. If the symmetric semigroup $\N$ is
generated by $4$ elements, say $\N =\ <\!\!m,m_1,m_2,m_3\!\!>$, then by
using Pinkham's equivariant deformation theory \cite{Pi74},
complete intersection theory and
a quasi-homogeneous version of Buchsbaum-Eisenbud's structure theorem
for Gorenstein ideals of codimension $3$ (see
\cite[p.\ 466]{BE77}), one can deduce that the affine monomial
curve $\Spec\ \k[\N] = \Spec\ \k[t^m,t^{m_1},t^{m_2},t^{m_3}]$
can be negatively smoothed without any obstructions (see \cite{Bu80},
\cite{W79} \cite[Satz 7.1]{W80}), hence
$\dim \M = \dim \proj(T^{1,-}_{\k[\N]|\k})$, and therefore
$$
\CM = \proj(T^{1,-}_{\k[\N]|\k}) 
$$
and so $\M$ is a dense open subvariety of $\proj(T^{1,-}_{\k[\N]|\k})$.
However, if $\N$ is generated by more than four elements, then
$\CM$ tends to be a proper subspace
of $\proj(T^{1,-}_{\k[\N]|\k})$, as documented in the above table
and discussed in the next section.


\section{Working with a family of symmetric semigroups}\label{6family}

\noindent In this section we apply our method to a family of symmetric
semigroups of multiplicity $6$ minimally generated by five elements.
For each positive integer $\tau$ we consider the semigroup
\begin{eqnarray*}
\N= & <6,2+6\tau,3+6\tau,4+6\tau,5+6\tau>\\ =& 6\Nat\,\sqcup\,
\displaystyle\bigsqcup_{j=2}^{5}(j+6\tau+6\Nat)\, \,\sqcup(7+12\tau+6\Nat) \ .
\end{eqnarray*}
Counting the numbers of gaps and picking up the largest nongap we obtain
$$g=6\tau+1 \ \ \text{and} \ \ \l_{g}=12\tau+1=2g-1$$
and so the semigroup is symmetric.

Let $\C$ be a complete integral Gorenstein curve, and $P$ a smooth point of $\C$
whose Weierstrass semigroup is equal to $\N$. As in Section \ref{construction}
we choose for each nongap $n\in\N$ a rational function $x_n$ on $\C$ with pole
divisor $n\,P$. We abbreviate
$$x:=x_6 \ \ \text{and} \ \  y_j:=x_{j+6\tau} \ \ (j=2,3,4,5)$$
and normalize
$$x_{6i}=x^i\ , \ \ x_{j+6\tau+6i}=x^iy_j \ \ \text{and} \ \
x_{7+12\tau+6i}=x^iy_2y_5$$
for each $i\geq 0$ and $j=2,3,4,5$. Now the $P$-hermitian basis
$\{x_{n_0},x_{n_1},\dots,x_{n_{g-1}}\}$ of the vector space $H^{0}(\C,(2g-2)P)$
of global sections of the canonical divisor $(2g-2)P=12\tau P$ consists of the
products
$$x^{0},\dots,x^{2\tau} \ \ \text{and} \ \ x^0y_j,\dots,x^{\tau-1}y_j
\ \ (j=2,3,4,5)\ .$$
Since $\l_2=2$, the complete integral Gorenstein curve $\C$ is nonhyperelliptic,
and so it can be identified with its image under the canonical embedding
$$(x_{n_0}:x_{n_1}\dots:x_{n_{g-1}}):\C\hookrightarrow \proj^{g-1}\ .$$
The projection map
$$(1:x:y_2:y_3:y_4:y_5):\C\rightarrow\proj^5$$
defines an isomorphism of the canonical curve $\C\subset\proj^{g-1}$ onto a
curve $\D\subset\proj^5$ of degree $6\tau+5$. The image $Q:=(0:0:0:0:0:1)\in\D$
of the distinguished Weier\-strass point $P:=(0:\dots:0:1)\in\C$ is the only point
of $\D$ that does not lie on the affine space $\AF^5\subset\proj^5$ of the
points with nonzero first coordinate.

To study quadratic relations of the canonical curve $\C\subset\proj^{g-1}$, we
consider the space of global sections of the bicanonical divisor
$(4g-4)P=24\tau P$. The $P$-hermitian basis $\{x_n\ | \ n\in\N$, $n\leq 4g-4\}$
of $H^0(\C,(4g-4)P)$ consists of the $3g-3$ functions
\begin{equation*}
\begin{array}{ll}
x^i & (i=0,1,\dots,4\tau) \\
x^iy_j & (i=0,1,\dots,3\tau-1\,, \ j=2,3,4,5) \\
x^iy_2y_5 & (i=0,1,\dots,2\tau-2)
\end{array}
\end{equation*}
which can be written as quadratic monomial expressions in
the projective coordinate functions
$x_{n_0},x_{n_1},\dots,x_{n_{g-1}}$. Let $X,Y_2,Y_3,Y_4$ and $Y_5$ be
indeterminates attached with the weights $6,2+6\tau,3+6\tau,4+6\tau$ and
$5+6\tau$, respectively. Having in mind the normalizations of the functions
$x_n$, we define for each nongap $n\in\N$ a monomial $Z_{n}$ of weight $n$ as
follows
$$Z_{6i}:=X^i\ , \ \ Z_{j+6\tau+6i}:=X^iY_j \ \ \text{and} \ \
Z_{7+12\tau+6i}:=X^iY_2Y_5\ .$$
Multiplying the functions $x,y_2,y_3,y_4,y_5$ by suitable constants, and
writing the nine products $y_iy_j$ with $(i,j)\neq(2,5)$ as a linear
combinations of the basis elements, we obtain nine polynomials in the
indeterminates $X,Y_2,Y_3,Y_4,Y_5$ that vanish identically on the affine curve
$\D\cap\AF^5=\D\setminus\{Q\}$, say
\begin{equation*}
\begin{array}{ll}
G_{i}=G_{i}^{(0)}-\displaystyle\sum_{j=1}^{12\tau+i}g_{ij}Z_{12\tau+i-j}
& (i=4,\dots,8) \\

F_{i}=F_{i}^{(0)}-\displaystyle\sum_{j=1}^{12\tau+i}f_{ij}Z_{12\tau+i-j}
& (i=6,8,9,10)
\end{array}
\end{equation*}
where
\begin{equation*}\label{9poly0}
 \begin{array}{lll}
G_{4}^{(0)}=Y_2^2-X^{\tau}Y_4 & G_5^{(0)}=Y_2Y_3-X^{\tau}Y_5 &
G_6^{(0)}=Y_3^2-X^{2\tau+1} \\

F_6^{(0)}=Y_2Y_4-X^{2\tau+1}
& G_7^{(0)}=Y_3Y_4-Y_2Y_5 &
G_8^{(0)}=Y_4^2-X^{\tau+1}Y_2 \\

F_8^{(0)}=Y_3Y_5-X^{\tau+1}Y_2
& F_9^{(0)}=Y_4Y_5-X^{\tau+1}Y_3 &
F_{10}^{(0)}=Y_5^2-X^{\tau+1}Y_4
 \end{array}
\end{equation*}
and where the summation index $j$ only varies through integers with
$12\tau+i-j\in\N$.

\begin{lem}
The ideal of the affine curve $\D\cap\AF^5$ is generated by
the nine polynomials $G_{i}$  $(i=4,\dots,8)$ and $F_{i}$ $(i=6,8,9,10)$.
\end{lem}
\begin{proof}
It follows by induction on descending degrees in $Y_2,\dots,Y_5$ that, modulo
the ideal generated by the nine polynomials, each polynomial in
$X,Y_2,\dots,Y_5$ is congruent to a polynomial whose terms are not divisible
by the nine products $Y_iY_j$ with $(i,j)\neq (2,5)$ i.e. which is a linear
combination of the monomials $Z_n$ of pairwise different weights $n\in\N$. Such
a linear combination $\sum c_nZ_n$ vanishes identically on the affine curve
$\D\cap\AF^5$ if and only if the corresponding linear combination
$\sum c_n x_{n}$ of the rational functions $x_n\in\k(\C)$ is equal to zero
i.e.\ $c_n=0$ for each $n\in\N$.
\end{proof}

If $\C$ is equal to the canonical monomial curve $\C^{(0)}\subset\proj^{g-1}$,
then the coefficients $g_{ij}$ and $f_{ij}$ are equal to zero, and so the ideal
of the affine monomial curve
$$
\D^{(0)}\cap\AF^5 =
\{(c^6,c^{2+6\tau},c^{3+6\tau},c^{4+6\tau},c^{5+6\tau})\ |\ c\in\k\}
$$
is generated by the nine quasi-homogeneous binomials $G_{i}^{(0)}$
$(i=4,\dots,8)$ and $F_{i}^{(0)}$ $(i=6,8,9,10)$.

In order to normalize some of the coefficients $g_{ij}$ and $f_{ij}$, we notice
that we have just the freedom to transform
\begin{equation*}
\begin{array}{lll}
x & \mapsto & c^6x+c_6 \\

y_2 & \mapsto & c^{2+6\tau}y_2+\sum_{i=0}^{\tau}c_{2+6i}x^{\tau-i} \\

y_3 & \mapsto & c^{3+6\tau}y_3+c_1y_2+\sum_{i=0}^{\tau}c_{3+6i}x^{\tau-i} \\

y_4 & \mapsto & c^{4+6\tau}y_4+c_1'y_3 + c_2'y_2 +
\sum_{i=0}^{\tau}c_{4+6i}x^{\tau-i} \\

y_5 & \mapsto & c^{5+6\tau}y_5+c_1''y_4 + c_2''y_3+c_3''y_2 +
\sum_{i=0}^{\tau}c_{5+6i}x^{\tau-i}
\end{array}
\end{equation*}
where $c\neq 0$, $c_j$, $c_j'$ and $c_j''$ are constants. We suppose that
the characteristic of the constant field is different from two and three.
Then we can normalize
$$f_{81}=g_{82}=f_{92}=f_{10,3}=0$$
(which are the only coefficients with $i-j\equiv 1\mod 6$),
$$g_{41}=g_{42}=g_{46}=0$$
and
$$f_{6,2+6i}=f_{8,3+6i}=f_{9,4+6i}=f_{9,5+6i}=0 \ \ (i=0,\dots,\tau) \ . $$
Now the isomorphism class of the pointed Gorenstein curve $(\C,P)$
determines uniquely the coefficients up to the following $\mathbb{G}_m$-action
$$g_{ij}\mapsto c^jg_{ij} \ \ \text{and} \ \ f_{ij}\mapsto c^jf_{ij} \ \
\text{where} \ \ c\in\mathbb{G}_m=\k^{*}\ .$$
We attach to the coefficients $g_{ij}$ and $f_{ij}$ the weight $j$. They
have to satisfy certain quasi-homogeneous polynomial equations,
which we will deduce from the syzygies of the affine curve
$\D\cap\AF^5$.

By applying the Syzygy Lemma we conclude that the five
quasi-homogeneous binomials
$$Z_{2g-2}G_{4}^{(0)},\ Z_{2g-2}G_{5}^{(0)},\
Z_{2g-2}(G_{6}^{(0)}-F_{6}^{0}),\ Z_{2g-2}G_{7}^{(0)} \ \text{and } \
Z_{2g-2}(G_{8}^{(0)}-F_{8}^{(0)})$$
of weight $2g-2+i+12\tau=24\tau+i$ where $i=4,5,6,7$ and $8$, respectively,
are linear combinations of the binomials $Z_nG_{j}^{(0)}$ $(j=4,\dots,8)$ and
$Z_nF_j^{(0)}$ $(j=6,8,9,10)$ with $n=2g-2+i-j<2g-2$ (and therefore
$X^{\tau-1}$ divides $Z_n$). More explicitly, we write up five syzygies of
the affine monomial curve $\D^{(0)}\cap\AF^5$:
\begin{equation*}
 \begin{array}{ll}
   X^{\tau+1}G_4^{(0)}-Y_4F_6^{(0)}+Y_2G_8^{(0)}=0 \\
   X^{\tau+1}G_{5}^{(0)}-Y_5F_{6}^{(0)}+Y_2F_9^{(0)}=0 \\
   X^{\tau+1}G_6^{(0)}-X^{\tau+1}F_{6}^{(0)}-Y_4F_8^{(0)}+Y_3F_9^{(0)}=0 \\
   X^{\tau+1}G_7^{(0)}-Y_5F_8^{(0)}+Y_3F_{10}^{(0)}=0 \\
   X^{\tau+1}G_8^{(0)}-X^{\tau+1}F_8^{(0)}-Y_5F_9^{(0)}+Y_4F_{10}^{(0)}=0
 \end{array}
\end{equation*}

\begin{rmk}\label{rem}
\textsl{
Actually, the Syzygy Lemma assures the existence of certain
$\frac{1}{2} (g-2)(g-5)$ syzygies of the canonical monomial curve
$\C^{(0)}\subset\proj^{g-1}$. However, using the equations $Z_{n+6}=X\cdot Z_n$
and factoring out powers of $X$ we can reduce to the five syzygies listed
above.
}
\end{rmk}
The five syzygies of the monomial curve $\D^{(0)}\cap\AF^5$ give rise to
five syzygies of the curve $\D\cap\AF^5$:

%







\medskip

\hspace{.14in}
$X^{\tau+1}G_4 - Y_4F_6 + Y_2G_8 \ =$ 

\noindent
 $\sum_{i=0}^{\tau}X^{\tau-i}
\left( \cc[1+6i]F_9 + \cc[3+6i]G_7 + (\cc[4+6i] {-} \kk[4+6i])F_6 -
\kk[5+6i]G_5 - \kk[6+6i]G_4 \right)$

\medskip

\hspace{.14in}
$X^{\tau+1}G_{5} - Y_5F_{6} + Y_2F_9 \ =$  

\noindent
   $\sum_{i=0}^{\tau} X^{\tau-i}\left( \cc[1+6i]F_{10} + \cc[3+6i]F_8
 - \m[6+6i]G_5 - \m[7+6i]G_4 \right)$

\medskip

\hspace{.14in}
  $X^{\tau+1}G_6 - X^{\tau+1}F_{6} - Y_4F_8 + Y_3F_9 \ =$ 

\noindent
  $\sum_{i=0}^{\tau}X^{\tau-i} ( \h[4+6i]G_8 + \h[5+6i]G_7
 - \m[6+6i]G_6 + \h[6+6i]F_6 - \m[7+6i]G_5)$

\medskip

\hspace{.14in}
 $X^{\tau+1}G_7 - Y_5F_8 + Y_3F_{10} \ =$

\noindent
  $\sum_{i=0}^{\tau}X^{\tau-i} ( \h[4+6i]F_9 {+} (\h[5+6i] {-} \n[5+6i])F_8
 {-}  \n[6+6i]G_7 {-} \n[7+6i]G_6 {-} \n[8+6i]G_5 )$

\medskip

\hspace{.14in}
  $X^{\tau+1}G_8-X^{\tau+1}F_8 - Y_5F_9 + Y_4F_{10} \ =$ 

\noindent
 $\sum_{i=0}^{\tau} X^{\tau-i} ( \m[6+6i]F_8
 - \n[5+6i]F_9 - \n[6+6i]G_8 - \n[7+6i]G_7 - \n[8+6i]F_6)$

\medskip

\noindent Indeed, by construction each right hand side differs from the
corresponding left hand side by a linear combination of the monomials
$Z_n$ that
vanishes identically on the curve $\D\cap\AF^5\cong\C\cap\AF^{g-1}$ and hence
is identically zero.

The vanishing of the coefficients of the five linear combinations provides us
with quasi-homogeneous equations between the coefficients $g_{ij}$ and
$f_{ij}$.
To express these equations in a concise manner, we introduce
polynomials in only one variable
$$g_i :=\sum_{r=1}^{12\tau+i}g_{ir}t^r =
G_{i}(t^{-6},t^{-2-6\tau},t^{-3-6\tau},t^{-4-6\tau},t^{-5-6\tau})t^{i+12\tau}
\ \
 (i=4,\dots,8)$$
and write each one as a sum of its \textit{partial polynomials}
$$\g[i,j]:=\sum_{r\equiv j \ \text{mod} \ 6}g_{ir}t^{r} \ \ (j=1,\dots,6)$$
which are defined by collecting terms whose exponents are in the same residue
class modulo $6$. In a similar way we define the polynomials $f_{i}$
$(i=6,8,9,10)$ and its partial polynomials $\f[i,j]$. Due to our normalizations
some of the partial polynomials are equal to zero, remaining only $41$ ones.
More precisely, we can write:
\begin{equation*}
 \begin{array}{ll}
g_4=\g[4,1]+\g[4,2]+\g[4,4]+\g[4,5]+\g[4,6] &
\ \ f_6=\f[6,1]+\f[6,3]+\f[6,4]+\f[6,6] \\
g_5=\g[5,1]+\g[5,2]+\g[5,3]+\g[5,5]+\g[5,6] &
\ \ f_8=\f[8,2]+\f[8,4]+\f[8,5]+\f[8,6]\\
g_6=\g[6,1]+\g[6,2]+\g[6,3]+\g[6,4]+\g[6,6] & \ \ f_9=\f[9,1]+\f[9,3]+\f[9,6] \\
g_7=\g[7,1]+\g[7,2]+\g[7,3]+\g[7,4]+\g[7,5] &
\ \ f_{10}=\f[10,1]+\f[10,2]+\f[10,4]+\f[10,5]+\f[10,6] \\
g_8=\g[8,2]+\g[8,3]+\g[8,4]+\g[8,5]+\g[8,6] &
\end{array}
\end{equation*}
The partial polynomials $\g[i,j]$ and $\f[i,j]$ with 
$i=j$ and $i=j+6$ \  i.e.
$\g[4,4]$, $\g[5,5]$, $\g[6,6]$, $\f[6,6]$ and  $\g[7,1]$, $\g[8,2]$, $\f[8,2]$,
$\f[9,3]$, $\f[10,4]$ have formal degree $i+12\tau$. The formal degree of
$\g[4,5]$, $\g[4,6]$, $\g[5,6]$, respectively $\f[9,1]$, $\f[10,1]$, $\f[10,2]$,
is equal to $j+6(\tau-1)$, respectively $j+6(\tau+1)$. The remaining $26$
partial polynomials have formal degree $j+6\tau$. Thus the number of the
coefficients that are still involved is equal to
$$4(2\tau+1)+5(2\tau+2)+3\tau+3(\tau+2)+26(\tau+1)-3=50\tau+43$$
where the discount by the number $3$ is due to the $3$ normalizations
$g_{41}=g_{42}=g_{46}=0$.
Now applying Theorem \ref{moduli} and Remark \ref{rem} we obtain an explicit
description of the compactified moduli space $\CM$.

\begin{thm}\label{modulifam}
Let $\N$ be a semigroup generated by $6$, $2+6\tau$, $3+6\tau$, $4+6\tau$ and
$5+6\tau$ where $\tau$ is a positive integer.
The isomorphism classes of the pointed complete
integral Gorenstein curves with Weierstrass semigroup $\N$ correspond
bijectively to the orbits of the $\mathbb{G}_m$-action on the quasi-cone of the
vectors of length $50\tau+43$ whose coordinates are the coefficients $g_{ij}$
and $f_{ij}$ of the $41$ partial polynomials that satisfy the five equations:

\begin{tabular}{l}
 $g_4-f_6+g_8=\f[6,1]f_9+\f[6,3]g_7+(\f[6,4]{-}\g[8,4])f_6-\g[8,5]g_5-\g[8,6]g_4$\\
 $g_5-f_6+f_9=\f[6,1]f_{10}+\f[6,3]f_8-\f[9,6]g_5-\f[9,1]g_4$ \\
 $g_6-f_6-f_8+f_9=\f[8,4]g_8+\f[8,5]g_7-\f[9,6]g_6+\f[8,6]f_6-\f[9,1]g_5$\\
 $g_7-f_8+f_{10}=\f[8,4]f_9+(\f[8,5]{-}\f[10,5])f_8-\f[10,6]g_7-\f[10,1]g_6-\f[10,2]g_5$\\
 $g_8-f_8-f_9+f_{10}=\f[9,6]f_8-\f[10,5]f_9-\f[10,6]g_8-\f[10,1]g_7-\f[10,2]f_6$
 \end{tabular}

\end{thm}
Thus the compactified moduli space $\CM$ admits an embedding into a weighted
projective space space of dimension $50\tau+42$.
To diminish the dimension of the ambient space, we project onto 
spaces of lower dimensions by eliminating some of the coordinates.

The $5$ equations of Theorem \ref{modulifam} can be rewritten in terms of $30$
polynomial equations between the $41$ partial polynomials. Among these
equations
there are $5$ linear ones, which we use to eliminate $5$ partial polynomials as
follows:
$$\f[6,4]=0, \ \ \g[8,3]=\f[6,3], \ \ \f[8,5]=0,\ \ \f[10,1]=\f[9,1],
\ \ \f[10,6]=\f[8,6]$$

\noindent There remain 25 inhomogeneous equations of degree two between
$36$ partial polynomials, among them the following $5$ equations:

\medskip
\begin{tabular}[h!]{l}
$\g[4, 2]+\g[8, 2] = \f[6, 1]\f[9, 1]+\f[6,3]\g[7,5]-\g[8, 5]\g[5, 3]-\g[8, 6]\g[4,2]$ \\

$\f[9, 3]-\f[6, 3] = \f[6, 1]\f[10, 2]+\f[6, 3]\f[8, 6]+\f[9, 1]\g[7, 2]$ \\

$\g[5, 6]-\f[6, 6]+\f[9, 6] =\f[6, 1]\f[10, 5]-\f[9,6]\g[5, 6]- \f[9, 1]\g[4, 5]$ \\

$\g[6, 2] -\f[8, 2]= \f[8, 4]\g[8, 4]-\f[9, 6]\g[6, 2]-\f[9, 1]\g[5, 1]$ \\

$\g[7, 4]-\f[8, 4]+\f[10, 4]=\f[8, 4]\f[9, 6]-\f[8, 6]\g[7, 4]-\f[9,1]\g[6, 3]
-\f[10, 2]\g[5, 2]$ 

\end{tabular}
\medskip

\noindent Since the formal degrees of $\g[8,2]$, $\f[9,3]$, $\f[6,6]$, $\f[8,2]$ 
and $\f[10,4]$ are equal to the formal degrees of the corresponding
equations, we can eliminate these five partial polynomials, remaining only $20$
equations between $31$ partial polynomials. As can be read off from the formal
degrees of the partial polynomials, this can be rephrased in terms of
$45\tau+40$ quasi-homogeneous equations between $35\tau+28$ coefficients. Some
of these equations may be identically zero.

Thus the compactified moduli space $\CM$ can be realized as a closed subvariety,
or more precisely, as an intersection of at most $45\tau+40$ hypersurfaces in
the weighted projective space of dimension $35\tau+27$, whose quasi-homogeneous
coordinates are the coefficients of the remaining $31$ partial polynomials.
We can continue in eliminating coefficients, until the remaining quasi-homogeneous
equations do not admit linear terms. However, if we would do this procedure
in an explicit way, our discussion would become very involved.

We first determine the weighted vector space $T^{1,-}_{\k[\N]|\k}$,
which is (up to an isomorphism) the locus of the linearizations of the
$30$ equations between the partial polynomials.  We can solve this system
of linear equations as follows:
\begin{small}
\begin{equation*}
 \begin{array}{lllll}
\f[6, 1] = \g[4, 1] & \f[9, 1] = \g[4, 1]-\g[5, 1] & \g[6, 1] = \g[5,
1] & \g[7, 1] = \g[5, 1]-\g[4, 1]
& \f[10,1]=\g[4, 1]-\g[5, 1]\\
\f[8, 2] = \g[6, 2] & \f[10, 2] = \g[4, 2]+\g[6, 2] & \g[5, 2] = 0 &
\g[7, 2] =-\g[4, 2] & \g[8, 2]=-\g[4, 2] \\
\f[9, 3] = \f[6, 3] & \g[5, 3] = 0 &  \g[6, 3] = 0 &  \g[7, 3] = 0 &
\g[8, 3] = \f[6, 3] \\
\f[6, 4] = 0 & \f[8, 4] = \g[6, 4] & \g[4, 4] = -\g[7, 4] & \g[8, 4] =
\g[7, 4] & \f[10, 4] = \g[6, 4]-\g[7, 4] \\
\f[8, 5] = 0 & \f[10, 5] = \g[4, 5] & \g[5, 5] = 0 & \g[7, 5] = -\g[4,
5] & \g[8, 5] = -\g[4, 5] \\
\f[10, 6]= \f[8, 6] & \f[6, 6] = \g[4, 6]{+}\g[8, 6] & \f[9, 6] =
\g[8, 6] & \g[5, 6] = \g[4, 6] & \g[6, 6] = \g[4, 6]{+}\f[8, 6]
 \end{array}
\end{equation*}\end{small}

\noindent Here we had to make choices, which partial polynomials should be
eliminated. However, we had to take care that the formal degrees on the left are
not smaller than the corresponding ones on the right. Now $T^{1,-}_{\k[\N]|\k}$
can be identified with the space of the vectors whose entries are the
coefficients of the remaining $11$ partial polynomials
$$\g[4,1],\ \g[5,1],\ \g[4,2],\ \g[6,2],\ \f[6,3],\ \g[6,4],\ \g[7,4],\
\g[4,5],\ \g[4,6],\ \g[8,6] \ \text{ and } \ \f[8,6]\ .$$
The only conditions the entries have to satisfy are the three normalizations
$$g_{41}=g_{42}=g_{46}=0\ .$$
Counting the coefficients that are still involved, we obtain
$$\dim T^{1,-}_{\k[\N|\k]}=11\tau+6 \ .$$
More precisely, counting the coefficients of a given weight $j$, we obtain the
dimension of the graded component of $T^{1}_{\k[\N]|\k}$ of negative weight
$-j$:

\begin{equation*}
\begin{array}{ll}
\dim T^{1}_{-1-6i}=\left\lbrace\begin{array}{ll} 1 & \text{ if } i=0
\\ 2 & (i=1,\dots,\tau) \end{array} \right. &
\dim T^{1}_{-2-6i}=\left\lbrace\begin{array}{ll} 1 & \text{ if } i=0
\\ 2 & (i=1,\dots,\tau) \end{array}\right. \\ & \\
\dim T^{1}_{-3-6i}= 1 \ \ (i=0,\dots,\tau) &
\dim T^{1}_{-4-6i}= 2 \ \ (i=0,\dots,\tau) \\ & \\
\dim T^{1}_{-5-6i}= 1 \ \ (i=0,\dots,\tau-1) &
\dim T^{1}_{-6-6i}=\left\lbrace\begin{array}{ll} 2 & \text{ if } i=0
\\ 3 & (i=1,\dots,\tau-1) \\ 2 & \text{ if }i=\tau \end{array}\right.
\end{array}
\end{equation*}
In the remaining cases the dimension of $T^{1}_{-j}$ is equal to zero.

Thus the compactified moduli space $\CM$ has been realized as a closed
subvariety of the $(11\tau+5)$-dimensional weighted projective space
$\proj(T^{1,-}_{\k[\N]|\k})$. It is cut out by $21\tau+28$ quasi-homogeneous
equations, which do not admit linear terms.

As discussed below, it will be much less expensive to obtain
the equations and the dimension of the quadratic quasi-cone $\SQ$
than the ones of the moduli variety $\CM$. In particular we can
make the eliminations in an explicit way.

To determine the quadratic quasi-cone $\SQ$, we just enter with our
solution of the system of $30$ linear equations into the quadratic terms
of the original $30$ equations of degree $\leq 2$, and eliminate the
same partial polynomials as in the linear case. We obtain only $5$
equations

\medskip
\begin{tabular}[h!]{l}
$\f[9, 1] = \g[4, 1]-\g[5, 1]+\g[4, 1](\f[8, 6]{-}\g[4, 6])
-\g[5,1](\g[8, 6]{-}\g[4, 6])+\f[6,3]\g[6, 4]$ \\

$\f[10, 2] = \g[4, 2]+\g[6, 2]+\g[6,2](\g[8, 6]{-}\g[4, 6])
+\g[4,2](\f[8,6]{-}\g[4, 6])-\g[6, 4]\g[7, 4]$ \\

$\g[6, 3] = -\g[4, 2]\g[5, 1]-\g[4, 1]\g[6, 2]-\g[6, 4]\g[4, 5]$ \\

$\g[8, 5] = -\g[4, 5]-\g[4, 5](\g[8, 6]{-}\g[4, 6])
-\g[7, 4]\g[4,1]-\f[6, 3]\g[4, 2]$ \\

$\f[10, 5] = \g[4, 5]+\g[4, 5](\f[8, 6]{-}\g[4, 6])
+\g[7, 4]\g[5,1]-\f[6, 3]\g[6, 2]$
\end{tabular}
\medskip

\noindent where the formal degrees on the left hand side are smaller than the
corresponding formal degrees of the right hand side, while in the
remaining $25$ equations the formal degrees on the left are
sufficiently large.
Thus the quadratic quasi-cone $\SQ$ is the subvariety of
$T^{1,-}_{\k[\N]|k}$ given by the five conditions

\begin{equation*}
\begin{array}{l}
\pi_{7+6\tau}\left(\g[4,1]\ff[8,6]-\g[5,1]\gg[8,6]+\f[6,3]\g[6,4]\right)=0 \\

\pi_{8+6\tau}\left(\g[6,2]\gg[8,6]+\g[4,2]\ff[8,6]-\g[6, 4]\g[7, 4]\right)=0 \\

\pi_{8+6\tau}\left(\g[4, 2]\g[5, 1]+\g[4, 1]\g[6, 2]+\g[6, 4]\g[4,
5]\right)=0 \\

\pi_{5+6\tau}\left(\g[4,5]\gg[8,6]+\g[7, 4]\g[4, 1]+\f[6, 3]\g[4, 2]\right)=0 \\

\pi_{5+6\tau}\left(\g[4,5]\ff[8,6]+\g[7, 4]\g[5, 1]-\f[6, 3]\g[6, 2]\right)=0
\end{array}
\end{equation*}
\medskip

\noindent where
$$ \ff[8,6]:=\f[8,6]-\g[4,6] \ \ \ \text{and} \ \ \ \gg[8,6]:=\g[8,6]-\g[4,6]$$
and where $\pi_i$ denotes the projection operator on the polynomials in $t$
that annihilates the terms of degree not larger than $i$.

We notice
that the five conditions do not depend on the $\tau+6$ coefficients $g_{51}$,
$g_{62}$, $f_{63}$, $g_{44}$, $g_{74}$, $\tilde{f_{86}}$, $\tilde{g_{86}}$ and
$g_{4,6i}$ $(i=2,\dots,\tau)$.
The five conditions on the remaining $10\tau$ coefficients can be expressed in
an elegant way in terms of five polynomial equations between ten elements of
the $\tau$-dimensional artinian algebra
\begin{equation*}
\A:=\k[\varepsilon]=\bigoplus\limits_{j=0}^{\tau-1}\,\k\,\varepsilon^{j} \ \
\text{where} \ \ \ \varepsilon^{\tau}=0\ .
\end{equation*}

\begin{thm}\label{aq6}
The quadratic quasi-cone $\SQ$ is isomorphic to the direct product
$$\SQ=V\times W$$
where $V$ is the $(\tau+6)$-dimensional weighted vector space with the weights
$1$, $2$, $3$, $4$, $4$, $6$ and $6i$ $(i=1,\dots,\tau)$ and where $W$ is the
quasi-cone consisting of the vectors
$$(\omega_1,\dots,\omega_{10})=
\left(\sum_{j=0}^{\tau-1}w_{1j}\,\varepsilon^j\,
,\dots,\,\sum_{j=0}^{\tau-1}w_{10,j}\,\varepsilon^j\right)
\in \A^{10}$$
satisfying the $5$ equations:
\begin{equation*}
\begin{array}{ll}
\omega_{1}\omega_{9}+\omega_{5}\omega_{6}-\omega_{2}\omega_{10} & = \ 0 \\

\omega_{4}\omega_{10}+\omega_{3}\omega_{9}-\omega_{6}\omega_{7} & = \ 0 \\

\omega_{1}\omega_{4}+\omega_{2}\omega_{3}+\omega_{6}\omega_{8} & = \ 0 \\

\omega_{8}\omega_{10}+\omega_{3}\omega_{5}+\omega_{1}\omega_{7} & = \ 0 \\

\omega_{2}\omega_{7}+\omega_{8}\omega_{9}-\omega_{4}\omega_{5} & = \ 0
\end{array}
\end{equation*}
in the artinian algebra $\A$. To the coefficients $w_{ij}$ are attached the
weights $\eta_i + 6(\tau-j)$ where $\eta_1,\dots,\eta_{10}=1,1,2,2,3,4,4,-1,6,6$.
\end{thm}
\begin{proof}
 We define
\begin{equation*}
\begin{array}{llll}
w_{1j}= g_{4,6\tau+1-6j}\, , &
w_{2j}= g_{5,6\tau+1-6j}\, , &
w_{3j}= g_{4,6\tau+2-6j}\, , &
w_{4j}= g_{6,6\tau+2-6j}\, , \\
w_{5j}= f_{6,6\tau+3-6j}\, , &
w_{6j}= g_{4,6\tau+4-6j}\, , &
w_{7j}= g_{7,6\tau+4-6j}\, , &
w_{8j}= g_{4,6\tau-1-6j}\, , \\

w_{9j}= \tilde{f}_{8,6\tau+6-6j}\, , &
w_{10,j}= \tilde{g}_{8,6\tau+6-6j}
\end{array}
\end{equation*}

\noindent
and notice that the five conditions on the $10\,\tau$ coefficients are equivalent
to the five quadratic equations in the artinian algebra $A$. \end{proof}

\begin{cor}\label{dimaq6}
$$\dim\SQ=8\,\tau+6\ .$$
\end{cor}
\begin{proof}
Since $\dim V=\tau+6$ we have to show that $$\dim W=7\,\tau\ .$$
For each $i=1,\dots,10$ let $W_i$ be the open subset of $W$ given by the
inequality $w_{i0}\neq0$, which means that $\omega_i$ is a unit in the local
artinian algebra $A$. If a vector $(\omega_1,\dots,\omega_{10})$ belongs to
$W_1$ then we can eliminate $\omega_9$, $\omega_4$ and $\omega_7$ from the
first, third and fourth quadratic equation, and the remaining two equations
become trivial. Thus $W_1$ has codimension $3\,\tau$ in $\A^{10}$ and
hence dimension $7\,\tau$. In a completely analogous way we see that
$$\dim W_i=7\,\tau \ \ (i=1,\dots,10)\ .$$
If $\tau=1$ then $W=W_1\cup\dots\cup W_{10}$ and therefore $\dim W=7$.

Now we assume that $\tau>1$. If a vector $(\omega_1,\dots,\omega_{10})\in W$
does not belong to the union $W_1\cup\dots\cup W_{10}$ i.e. $w_{ij}=0$
whenever $j=0$, then the ten coefficients $w_{ij}$ with $j=\tau-1$ do not enter
into the five quadratic equations, and by induction we obtain
$$\dim(W\setminus(W_1\cup\dots\cup W_{10}))=7(\tau-2)+10=7\,\tau-4<7\,\tau$$
and therefore $\dim W=7\,\tau$. \end{proof}

Now applying Theorem \ref{ubound} we obtain an upper bound for the dimension of
the moduli variety $$\dim\CM\leq 8\,\tau+5$$ which for each $\tau > 1$
is better than Deligne's bound $2\,g-1=12\,\tau+1$.

\end{document}